\documentclass[smallextended,final,envcountsect]{svjour3}

\usepackage[top=1in, bottom=1in, left=1.2in, right=1.2in]{geometry}

\usepackage{cite, amsmath, url, amssymb}
\usepackage[misc,geometry]{ifsym}
\usepackage[title]{appendix}
\usepackage[hang,small,bf]{caption}

\RequirePackage{fix-cm}
\smartqed
\journalname{JOTA}

\newcommand{\eg}{{\it e.g.}}
\newcommand{\ie}{{\it i.e.}}

\newcommand{\BA}{\begin{array}}
\newcommand{\EA}{\end{array}}
\newcommand{\K}{\mathcal{K}}

\newcommand{\ones}{\mathbf 1}
\newcommand{\reals}{{\mathbb{R}}} 

\newcommand{\Tr}{\mathop{\bf Tr}}
\newcommand{\diag}{\mathop{\bf diag}}

\newcommand{\argmin}{\mathop{\rm argmin}}

\newcommand{\minimize}{minimize}
\newcommand{\maximize}{maximize}
\newcommand{\subjectto}{s.t.~}
\newcommand{\setsep}{:}
\renewcommand{\S}{Sect.~}

\newcommand{\cones}{\mathcal{K}}

\newcommand{\e}[2]{$#1\times 10^{#2}$}
\newcommand{\eb}[2]{$\mathbf{#1\times 10^{#2}}$}

\title{Conic Optimization via Operator Splitting and
Homogeneous Self-Dual Embedding}

\author{Brendan O'Donoghue \and Eric Chu \and Neal Parikh \and Stephen Boyd}
\date{\today}

\institute{Brendan O'Donoghue \Letter\ \email {bodonoghue85@gmail.com},
Department of Electrical Engineering, Stanford University \and
Eric Chu, Department of Electrical Engineering, Stanford University \and
Neal Parikh, Department of Computer Science, Stanford University \and
Stephen Boyd, Department of Electrical Engineering, Stanford University}

\begin{document}

\maketitle

\begin{abstract}

We introduce a first order method for solving very large convex cone programs. 
The method uses an operator splitting method, the alternating directions method
of multipliers, to solve the homogeneous self-dual embedding, an equivalent
feasibility problem involving finding a nonzero point in the intersection of a
subspace and a cone.

This approach has several favorable properties. Compared to interior-point
methods, first-order methods scale to very large problems, at the
cost of requiring more time to reach very high accuracy.  Compared to other
first-order methods for cone programs, our approach finds both primal and dual
solutions when available or a certificate of infeasibility or unboundedness
otherwise, is parameter-free, and the per-iteration cost of the method is the
same as applying a splitting method to the primal or dual alone.

We discuss efficient implementation of the method in detail, including direct
and indirect methods for computing projection onto the subspace, scaling the
original problem data, and stopping criteria.  We describe an
open-source implementation, 
which handles the usual (symmetric)
non-negative, second-order, and semidefinite cones as well
as the (non-self-dual) exponential and power cones and their duals.
We report numerical results that show speedups over
interior-point cone solvers for large problems, and scaling to
very large general cone programs.
\end{abstract}

\section{Introduction}
\label{s-intro}

In this paper we develop a method for solving convex cone
optimization problems that can (a) provide primal or dual certificates of
infeasibility when relevant and (b) scale to large problem sizes. The general
idea is to use a first-order method to solve the homogeneous self-dual
embedding of the primal-dual pair; the homogeneous self-dual embedding provides
the necessary certificates, and first-order methods scale well to large problem
sizes.

The homogeneous self-dual embedding is a single convex feasibility problem that
encodes the primal-dual pair of optimization problems. Solving the embedded
problem involves finding a nonzero point in the intersection of two convex
sets, a convex cone and a subspace.  If the original pair is solvable, then a
solution can be recovered from any nonzero solution to the embedding;
otherwise, a certificate of infeasibility is generated that proves that the
primal or dual is infeasible (and the other one unbounded).  
The homogeneous self-dual embedding has been
widely used with interior-point methods \cite{Ye:11,sedumi, SY:12}.

We solve the embedded problem with an operator splitting method known as the
\emph{alternating direction method of multipliers} (ADMM)
\cite{GlM:75,GaM:76,G:83,Eck:89}; see \cite{BP:11} for a recent survey.  It can
be viewed as a simple variation of the classical alternating projections
algorithm for finding a point in the intersection of two convex sets.
Roughly speaking, ADMM adds a dual state variable to the basic method, 
which can substantially improve convergence.
The overall method can
reliably provide solutions to modest accuracy after a relatively small number
of iterations and can solve large problems far more quickly than interior-point
methods. 
(It may not be suitable if high accuracy is required, due to the slow
`tail convergence' of first order methods in general, and ADMM in
particular \cite{HY:12b}.) To the best of our knowledge, this
is the first application of a first-order method to solving such embeddings.
The approach described in this paper combines a number of different ideas that
are well-established in the literature, such as cone programming and operator
splitting methods. We highlight various dimensions along which our method can
be compared to others.

Some methods for solving cone programs only return primal solutions, while
others can return primal-dual pairs. In addition, some methods can only handle
feasible problems, while other methods can also return certificates of
infeasibility or unboundedness. The idea of homogeneous self-dual embedding is
due to Ye and others \cite{YT:94,XH:96}.  Self-dual embeddings have generally
been solved via interior-point methods \cite{NeN:94}, while the literature on
other algorithms has generally yielded methods that cannot return certificates
of infeasibility; see, \eg, \cite{WGY:10,LLR:11,AI:13}.

Our approach involves converting a primal-dual pair into a convex feasibility
problem involving finding a point in the intersection of two convex sets. There
are many projection algorithms that could be used to solve this
kind of problem, such as the classical alternating directions method or
Dykstra's alternating projections method \cite{BD:86, BB:94}, amongst others
\cite{CCCH:12, CE:94}. For a further discussion of these and many other
projection methods, see Bauschke and Koch \cite{BK:13}.
Any of these methods could be used to solve the problem in 
homogeneous self-dual embedding form.

Operator splitting techniques go back to the 1950s; ADMM itself was developed
in the mid-1970s \cite{GlM:75,GaM:76}.  Since then a rich literature has
developed around ADMM and related methods \cite{G:83,Eck:89,EB:92,CP:12, C:13,
KP:14, GT:89, LM:79, G:84, FG:83}. 
Many equivalences exist between ADMM and other operator splitting methods.
It was shown in \cite{G:83} that ADMM is
equivalent to the variant of Douglas-Rachford splitting presented in
\cite{LM:79} (the original, more restrictive, form of Douglas-Rachford
splitting was presented in \cite{DR:56}) applied to the dual problem, which
itself is equivalent to Rockafellar's proximal point algorithm \cite{R:76,
EB:92}.  

Douglas-Rachford splitting is also equivalent to Spingarn's `method of partial
inverses' \cite{Spi:83,Spi:85a, Spi:85b} when one of the operators is the
normal cone map of a linear subspace \cite{Eck:88, Eck:89}. In this paper we
apply ADMM to a problem where one of the functions is the indicator of a linear
subspace, so our algorithm can also be viewed as an application of Spingarn's
method. Another closely related technique is the `split-feasibility problem', which
seeks two points related by a linear mapping, each of which is constrained to
be in a convex set \cite{CE:94, B:02, CMS:07, C:01}.

In \cite{Eck:89} and \cite{MW:14} it was shown that equivalences exist between
ADMM applied to the primal problem, the dual problem, and a saddle point
formulation of the problem; in other words, ADMM is (in a sense) itself 
self-dual.

These techniques have been used in a broad range of applications including
imaging \cite{Com:96, GO:09, OV:13}, control \cite{jovanovic:10, annergren:12,
oper_splt_ctrl,MX:12, BOD:12}, estimation \cite{admm_tv_est}, signal processing
\cite{CW:06,CoP:07,CP:09, YZ:11}, finance \cite{port_opt_bound}, distributed
optimization \cite{PB:12,K:13}, and many others.

There are several different ways to apply ADMM to solve cone programs
\cite{WGY:10, BP:11}.  In some cases, these are applied to the original cone
program (or its dual) and yield methods that can return primal-dual pairs,  but
cannot handle infeasible or unbounded problems.


The indirect version of our method interacts with the data solely by
multiplication by the data matrix or its adjoint, which we can informally refer
to as a `scientific computing' style algorithm; it is also called a
`matrix-free method'.  There are several other methods that share similar
characteristics, such as \cite{CP:11,BCG:10,Gon:12, MOS:13, MOS:14a, MOS:14b,
ZST:10, OC:12}, as well as some techniques for solving the split-feasibility
problem \cite{B:02}.  See Esser et al. \cite{EZC:10} for a detailed discussion
of various first-order methods and the relationships between them, and Parikh
and Boyd \cite{PB:14} for a survey of proximal algorithms in particular.

\paragraph{Outline.}
In \S\ref{s-conic-opt} we review convex cone
optimization, conditions for optimality, and the homogeneous self-dual
embedding.  In \S\ref{s-algorithm}, we derive an algorithm that solves
(\ref{e-coneprob}) using ADMM applied to the homogeneous self-dual embedding of
a cone program. In \S\ref{s-fast-impl}, we discuss how to perform the sub-steps
of the procedure efficiently.  In \S\ref{s-scaling} we introduce a scaling
procedure that greatly improves convergence in practice. We conclude
with some numerical examples in \S\ref{s-num-ex}, including (when applicable,
\ie, the problems are small enough and involve only symmetric cones) a
comparison of our approach with a state-of-the-art interior-point method, both
in quality of solution and solution time.

\section{Conic Optimization}
\label{s-conic-opt}
Consider the \emph{primal-dual pair} of (convex) cone optimization problems
\begin{equation}
\label{e-coneprob}
\BA{lcl}
\BA{l}
\BA{ll}
\mbox{\minimize} & c^T x
\EA \\
\BA{ll}
\mbox{\subjectto} & Ax + s = b \\
& (x,s) \in \reals^n \times \K,
\EA
\EA
& &
\BA{l}
\BA{ll}
\mbox{\maximize} & - b^T y
\EA \\
\BA{ll}
\mbox{\subjectto} & -A^Ty + r = c \\
& (r,y) \in \{0\}^n \times \K^*.
\EA
\EA
\EA
\end{equation}
Here $x \in \reals^n$ and $s \in \reals^m$ (with $n \leq m$) are the primal
variables, and $r \in \reals^n$ and $y \in \reals^m$ are the dual variables 
We refer to $x$ as the primal variable, $s$ as the primal slack variable, $y$
as the dual variable, and $r$ as the dual residual.  The set $\mathcal K$ is a
nonempty, closed, convex cone with dual cone $\mathcal K^*$, and $\{0\}^n$ is
the dual cone of $\reals^n$, so the cones $\reals^n \times \mathcal K$ and
$\{0\}^n \times \K^*$ are duals of each other.  The problem data are $A\in
\reals^{m \times n}$, $b\in \reals^m$, $c\in\reals^n$, and the cone $\K$.
(We consider all vectors to be column vectors.)

The primal and dual optimal values are denoted $p^\star$ and $d^\star$,
respectively; we allow the cases when these are infinite: $p^\star = +\infty$
($-\infty$) indicates primal infeasibility (unboundedness), and $d^\star =
-\infty$ ($+\infty$) indicates dual infeasibility (unboundedness).  It is easy
to show weak duality, \ie, $d^\star \leq p^\star$, with no assumptions on the
data.  We will assume that strong duality holds, \ie, $p^\star = d^\star$,
including the cases when they are infinite.

\subsection{Optimality Conditions}
\label{s-primal-dual}

When strong duality holds, the KKT (Karush-Kuhn-Tucker) conditions are
necessary and sufficient for optimality. Explicitly,
$(x^\star,s^\star,r^\star,y^\star)$ satisfies the KKT conditions, and so is
primal-dual optimal, when
\[
Ax^\star + s^\star = b, \quad s^\star \in \K, \quad
A^T y^\star + c = r^\star,\quad r^\star = 0, \quad y^\star \in \K^*, \quad
(y^\star)^T s^\star = 0,
\]
\ie, when $(x^\star,s^\star)$ is primal feasible, $(r^\star,y^\star)$ is dual
feasible, and the complementary slackness condition $(y^\star)^T s^\star = 0$
holds.  
The complementary slackness condition can equivalently be replaced by the condition
\[
c^T x^\star + b^T y^\star = 0,
\]
which explicitly forces the \emph{duality gap}, $c^T x + b^T y$, 
to be zero.

\subsection{Certificates of Infeasibility}

If strong duality holds, then exactly one of the sets
\begin{align}
\mathcal P &= \{ (x,s) \setsep Ax + s = b,\ s \in \K \} \label{e-primal-p-alt}, \\
\mathcal D &= \{ y \setsep A^T y = 0,\ y \in \K^*,\ b^T y < 0 \} \label{e-primal-d-alt},
\end{align}
is nonempty, a result known as a \emph{theorem of strong alternatives}
\cite[\S5.8]{BV:04}.  Since the set $\mathcal P$ encodes primal feasibility,
this implies that any dual variable $y \in \mathcal D$ serves as a \emph{proof}
or \emph{certificate} that the set $\mathcal P$ is empty, \ie, that the problem
is primal infeasible.  Intuitively, the set $\mathcal D$ encodes the
requirements for the dual problem to be feasible but unbounded.

Similarly, exactly one of the following two sets is nonempty:
\begin{align}
\tilde{\mathcal P} &= \{ x \setsep -Ax \in \K,\ c^T x < 0 \} \label{e-dual-p-alt}, \\
\tilde{\mathcal D} &= \{ y \setsep A^T y = -c,\ y \in \K^* \} \label{e-dual-d-alt}.
\end{align}
Any primal variable $x \in \tilde{\mathcal P}$ is a certificate of dual infeasibility.

\subsection{Homogeneous Self-Dual Embedding}
\label{s-extembed}

The original pair of problems~\eqref{e-coneprob} can be converted into
a single feasibility problem by embedding the KKT conditions into
a single system of equations and inclusions that the primal and dual
optimal points must jointly satisfy. The embedding is as follows:
\begin{equation}
\label{e-pdem}
    \begin{bmatrix} r \\ s \\ 0 \end{bmatrix} =
    \begin{bmatrix}
        0   & A^T \\
        -A  & 0 \\
        c^T & b^T
    \end{bmatrix}
    \begin{bmatrix} x \\ y \end{bmatrix} +
    \begin{bmatrix} c \\ b \\ 0 \end{bmatrix},
        \quad
    (x,s,r,y) \in \reals^n \times \K \times \{0\}^n \times \K^*.
\end{equation}
Any $(x^\star,s^\star,r^\star,y^\star)$ that satisfies \eqref{e-pdem} is
optimal for \eqref{e-coneprob}. However, if \eqref{e-coneprob} is primal or
dual infeasible, then \eqref{e-pdem} has no solution.

The homogeneous self-dual embedding \cite{YT:94}
addresses this shortcoming:
\begin{equation}
\label{e-hsd}
    \begin{bmatrix} r \\ s \\ \kappa \end{bmatrix} =
    \begin{bmatrix}
        0    & A^T  & c \\
        -A   & 0    & b \\
        -c^T & -b^T & 0
    \end{bmatrix}
    \begin{bmatrix} x \\ y \\ \tau \end{bmatrix},
        \quad
    (x,s,r,y,\tau,\kappa) \in \reals^n \times \K \times \{0\}^n \times \K^* \times \reals_+ \times \reals_+.
\end{equation}
This embedding introduces two new variables, $\tau$ and $\kappa$, that are
non-negative and complementary, \ie, at most one is nonzero.  To see
complementarity note that the inner product between $(x,y,\tau)$ and
$(r,s,\kappa)$ at any solution must be zero due to the skew symmetry of the
matrix in (\ref{e-hsd}), and the individual components $x^Tr$, $y^Ts$, and
$\tau \kappa$ must each be non-negative by the definition of dual cones.

The reason for using this embedding is that the different
possible values of $\tau$ and $\kappa$ encode the different possible outcomes.
If $\tau$ is nonzero at the solution, then it serves as a scaling factor that
can be used to recover the solutions to (\ref{e-coneprob}); otherwise, if
$\kappa$ is nonzero, then the original problem is primal or dual infeasible. In
particular, if $\tau = 1$ and $\kappa = 0$ then the self-dual embedding reduces to
the simpler embedding (\ref{e-pdem}). 

Any solution of the self-dual embedding $(x,s,r,y,\tau,\kappa)$ falls into one of three cases: 
\begin{enumerate}
\item $\tau > 0$ and $\kappa = 0$. The point 
\[
(\hat x, \hat y, \hat s) = (x/\tau, y/\tau, s/\tau)
\]
satisfies the KKT conditions of (\ref{e-coneprob}) and so is a
primal-dual solution.

\item $\tau = 0$ and $\kappa > 0$. This implies 
that the gap $c^T x + b^T y$ is negative, which
immediately tells us that the problem is either primal or dual infeasible.
\begin{itemize}
\item If $b^T y < 0$, then $\hat y = y/(b^T y)$
is a certificate of primal infeasibility (\ie, $\mathcal D$ is nonempty) since
\[
A^T \hat y = 0,\quad \hat y \in \mathcal{K}^*,\quad b^T \hat y = -1.
\]
\item If $c^T x < 0$, then $\hat x = x/(-c^T x)$
is a certificate of dual infeasibility (\ie, $\tilde{\mathcal P}$ is nonempty)
since
\[
-A\hat x \in \mathcal{K},\quad c^T\hat x = -1.
\]
\item If both $c^T x < 0$ and $b^T y < 0$, then
the problem is both primal and dual infeasible (but the strong duality
assumption is violated).
\end{itemize}
\item $\tau = \kappa = 0$. If one of $c^T x$ or $b^T y$
is negative, then it can be used to derive a certificate of primal or
dual infeasibility.
Otherwise nothing can be concluded about the original problem.
Note that zero is always a solution to (\ref{e-hsd}), but steps can be taken to
avoid it, as we discuss in Section \ref{s-convergence}.
\end{enumerate}

The system \eqref{e-hsd} is homogeneous because if
$(x,s,r,y,\tau,\kappa)$ is a solution to the embedding, then so is
$(tx,ts,tr,ty,t\tau,t\kappa)$ for any $t \geq 0$, and when 
$t>0$ this scaled value yields the same primal-dual solution or
certificates for (\ref{e-coneprob}). The embedding is also self-dual, which we
show below.

\paragraph{Notation.}
To simplify the subsequent discussion, let
\[
u = \begin{bmatrix} x\\y\\\tau \end{bmatrix}, \quad
v = \begin{bmatrix} r\\s\\\kappa \end{bmatrix}, \quad
Q = 
\begin{bmatrix}
0 & A^T & c \\
-A & 0 & b \\
-c^T & -b^T & 0
\end{bmatrix}.
\]
The homogeneous self-dual embedding \eqref{e-hsd} can then be expressed as
\begin{equation}
\label{e-hsdnice}
\begin{array}{ll}
\mbox{find} & (u,v)\\
\mbox{\subjectto} & v = Qu \\
& (u,v) \in \mathcal{C} \times \mathcal{C}^*,
\end{array}
\end{equation}
where $\mathcal{C} = \reals^n \times \cones^* \times \reals_+$ is a cone with
dual cone $\mathcal{C}^* = \{ 0 \}^n \times \cones \times \reals_+$. 
We are interested in finding a nonzero solution of
the homogeneous self-dual embedding (\ref{e-hsdnice}). 
In the sequel, $u_x, u_y, u_\tau$ and
$v_r,v_s, v_\kappa$ will denote the entries of $u$ and $v$ that correspond
to $x,y,\tau$ and $r,s,\kappa$, respectively.

\paragraph{Self-dual property.}
Let us show that the feasibility problem (\ref{e-hsdnice}) is self-dual.
The Lagrangian has the form
\[
L(u,v,\nu,\lambda,\mu) = \nu^T(Qu-v) - \lambda^T u -\mu ^T v,
\]
where the dual variables are $\nu,\lambda, \mu$, with
$\lambda \in \mathcal C^*$, $\mu \in \mathcal C$.
Minimizing over the primal variables $u,v$, we conclude that
\[
Q^T \nu - \lambda = 0, \quad -\nu - \mu = 0.
\]
Eliminating $\nu = -\mu$ and using $Q^T = -Q$ we can write the dual problem as
\[
\begin{array}{ll}
\mbox{find} & (\mu,\lambda)\\
\mbox{\subjectto} & \lambda = Q\mu \\
& (\mu,\lambda) \in \mathcal{C} \times \mathcal{C}^*,
\end{array}
\]
with variables $\mu,\lambda$.
This is identical to (\ref{e-hsdnice}).

\section{Operator Splitting Method}
\label{s-algorithm}

The convex feasibility problem~(\ref{e-hsdnice}) can be solved by many methods,
ranging from simple alternating projections to sophisticated interior-point
methods.  We are interested in methods that scale to very large problems, so we
will use an operator splitting method, the alternating direction method of
multipliers (ADMM).  There are many operator splitting methods (some of which
are equivalent to ADMM) that could
be used to solve the convex feasibility problem, such as Douglas-Rachford
iteration, split feasibility methods, Spingarn's method of partial inverses,
Dykstra's method, and others.  While we have not tried these other methods, we
suspect that many of them would yield comparable results to ADMM.  Moreover,
much of our discussion below, on simplifying the iterations and efficiently
carrying out the required steps, would also apply to (some) other operator
splitting methods.

\subsection{Basic Method}
ADMM is an operator splitting method that can solve convex problems
of the form
\begin{equation}
\label{e-admmprob}
\begin{array}{ccccc}
\mbox{\minimize} & \left[f(x) + g(z)\right] & & \mbox{\subjectto} & x = z.
\end{array}
\end{equation}
(ADMM can also solve problems where $x$ and $z$ are affinely related; 
see \cite{BP:11} and the references therein.)
Here, $f$ and $g$ may be nonsmooth or take on infinite values to encode
implicit constraints. 
The basic ADMM algorithm is
\begin{align*}
x^{k+1} &= \argmin_x \left(f(x) + (\rho/2) \|x - z^k - \lambda^k \|_2^2\right) \\
z^{k+1} &= \argmin_z \left(g(z) + (\rho/2) \|x^{k+1} - z - \lambda^k \|_2^2\right) \\
\lambda^{k+1} &= \lambda^k - x^{k+1} + z^{k+1},
\end{align*}
where $\rho>0$ is a step size parameter and $\lambda$ is the (scaled) dual
variable associated with the constraint $x=z$, and the superscript $k$ denotes
iteration number.  The initial points $z^0$ and $\lambda^0$ are arbitrary, but
are usually taken to be zero.  Under some very mild conditions
\cite[\S3.2]{BP:11}, ADMM converges to a solution, in the following sense:
$f(x^k) + g(z^k)$ converges to the optimal value, $\lambda^k$ converges to an
optimal dual variable, and $x^k - z^k$, the equality constraint residual,
converges to zero. Additionally, for the restricted form we consider in
(\ref{e-admmprob}), we have the stronger guarantee that $x^k$ and $z^k$ converge
to a common value; see, \eg, \cite[\S 5]{EB:92}.  We will mention later some
variations on this basic ADMM algorithm with similar convergence guarantees.

To apply ADMM, we transform the embedding \eqref{e-hsdnice} to ADMM form
\eqref{e-admmprob}:
\begin{equation}
\label{e-consensus}
\begin{array}{ccccc}
\mbox{\minimize} & \left[I_{\mathcal{C} \times \mathcal{C}^*}(u, v) + I_{Q u= v}(\tilde u,\tilde v)\right]
& & \mbox{s.t.~} & (u,v) = (\tilde u, \tilde v),
\end{array}
\end{equation}
where $I_\mathcal{S}$ denotes the indicator function
\cite[\S4]{Roc:70} of the set $\mathcal{S}$. A direct application of ADMM to
the self-dual embedding, written as (\ref{e-consensus}), yields the following algorithm:
\begin{equation}
\begin{array}{rcl}
(\tilde u^{k+1}, \tilde v^{k+1}) &=& \Pi_{Qu = v}(u^k + \lambda^k, v^k + \mu^k) \\[1ex]
u^{k+1} &=& \Pi_\mathcal{C}(\tilde u^{k+1} - \lambda^k) \\[1ex]
v^{k+1} &=& \Pi_\mathcal{C^*}(\tilde v^{k+1} - \mu^k) \\[1ex]
\lambda^{k+1} &=& \lambda^k -\tilde  u^{k+1} + u^{k+1} \\[1ex]
\mu^{k+1} &=& \mu^k - \tilde v^{k+1} + v^{k+1},
\end{array}
\label{e-SCS-naive}
\end{equation}
where $\Pi_{\mathcal S}(x)$ denotes the Euclidean projection of $x$ onto the
set $\mathcal S$. Here, $\lambda$ and $\mu$ are dual variables for
the equality constraints on $u$ and $v$, respectively.

\subsection{Simplified Method}

In this section we show that the basic ADMM algorithm (\ref{e-SCS-naive}) 
given above can be simplified using properties of our specific problem.

\subsubsection{Eliminating Dual Variables}

If we initialize $\lambda^0 = v^0$ and $\mu^0 = u^0$, then $\lambda^k = v^k$
and $\mu^k = u^k$ for all subsequent iterations. This result allows us to
eliminate the dual variable sequences above. This will also simplify the linear
system in the first step and remove one of the cone projections.

\begin{proof}
The proof is by induction. The base case holds because we can initialize the
variables accordingly. Assuming that $\lambda^k = v^k$ and $\mu^k = u^k$, the
first step of the algorithm becomes
\begin{equation}
\label{e-projA}
(\tilde u^{k+1}, \tilde v^{k+1}) = \Pi_\mathcal{Q}\left(u^k + \lambda^k, v^k + \mu^k\right) 
= \Pi_\mathcal{Q}\left( u^k + v^k, u^k + v^k\right),
\end{equation}
where $\mathcal{Q} = \{(u,v) \setsep Qu = v \}$. 

The orthogonal complement of $\mathcal Q$ is $\mathcal{Q^\perp} = \{(v,u) \setsep
Qu = v \}$ because $Q$ is skew-symmetric. 
It follows that if $(u,v) = \Pi_\mathcal{Q}(z,z)$, then $(v,u) =
\Pi_\mathcal{Q^\perp}(z,z)$ for any $z$, since the two projection problems are
identical save for reversed output arguments. This implies that
\begin{equation}
\label{e-projAp}
(\tilde v^{k+1}, \tilde u^{k+1}) = \Pi_\mathcal{Q^\perp}\left(u^k + v^k, u^k + v^k\right).
\end{equation}
Recall that $z = \Pi_\mathcal{Q}(z) + \Pi_\mathcal{Q^\perp}(z)$ for any $z$.
With (\ref{e-projA}) and (\ref{e-projAp}), this gives
\begin{equation}
\label{e-dsum}
u^k +v^k = \tilde u^{k+1} + \tilde v^{k+1}.
\end{equation}

The \emph{Moreau decomposition} \cite[\S2.5]{PB:14} of $x$ with
respect to a nonempty, closed, convex cone $\mathcal C$ is given by
\begin{equation}\label{e-moreau}
x = \Pi_\mathcal{C}(x) + \Pi_\mathcal{-C^*}(x),
\end{equation}
and moreover, the two terms on the right-hand side are orthogonal.
It can be written equivalently as
$x = \Pi_\mathcal{C}(x) - \Pi_\mathcal{C^*}(-x)$.
Combining this with \eqref{e-dsum} gives
\begin{align*}
u^{k+1} &= \Pi_\mathcal{C}(\tilde u^{k+1} - v^k) \\
&= \Pi_\mathcal{C}(u^k - \tilde v^{k+1}) \\
&= u^k - \tilde v^{k+1} + \Pi_{\mathcal{C}^*}(\tilde v^{k+1} - u^k) \\
&= u^k - \tilde v^{k+1} + v^{k+1} \\
&= \mu^{k+1}.
\end{align*}
A similar derivation yields $\lambda^{k+1}= v^{k+1}$, which completes
the proof. This lets us eliminate the sequences $\lambda^k$ and $\mu^k$. \qed
\end{proof}

Once the value $u^{k+1} =
\Pi_\mathcal{C}(\tilde u^{k+1} - v^k)$ has been calculated, the step that projects
onto the dual cone $\mathcal{C}^*$ can be replaced with
\[
v^{k+1} =  v^k - \tilde u^{k+1} +  u^{k+1}.
\]
This follows from the $\lambda^k$ update, which is typically cheaper than a
projection step. Now no sequence depends on $\tilde v^k$ any longer, so it
too can be eliminated.

\subsubsection{Projection Onto Affine Set}

Each iteration, the algorithm \eqref{e-SCS-naive} computes a
projection onto $\mathcal{Q}$ by solving
\[
\begin{array}{ccccc}
\mbox{\minimize} & \left[(1/2) \|u - u^k - v^k\|_2^2 + (1/2)\|v - u^k - v^k\|_2^2\right] &&  
\mbox{\subjectto} & v = Qu
\end{array}
\]
with variables $u$ and $v$.
The KKT conditions for this problem are
\begin{equation}
\begin{bmatrix}
I & Q^T \\ Q & -I
\end{bmatrix}
\begin{bmatrix}
u \\ \mu
\end{bmatrix}
=
\begin{bmatrix}
u^k + v^k \\ u^k + v^k
\end{bmatrix},
\label{e-kkt}
\end{equation}
where $\mu \in \reals^{m+n+1}$ is the dual variable associated with the
equality constraint $Qu - v = 0$. 
By eliminating $\mu$, we obtain 
\[
\tilde u^{k+1} = (I + Q^T Q)^{-1}(I-Q)(u^k + v^k). 
\]
The matrix $Q$ is skew-symmetric, so this simplifies to
\[
\tilde u^{k+1} = (I+Q)^{-1}(u^k + v^k).
\]
(The matrix $I + Q$ is guaranteed to be invertible since $Q$ is 
skew-symmetric.)

\subsubsection{Final Algorithm}

Combining the simplifications of the previous sections,
the final algorithm is
\begin{equation}
\begin{array}{rcl}
\tilde u^{k+1} &=& (I + Q)^{-1} (u^k + v^k )\\[1ex]
u^{k+1} &=& \Pi_\mathcal{C}\left(\tilde u^{k+1} - v^k\right) \\[1ex]
v^{k+1} &=&  v^k - \tilde u^{k+1} + u^{k+1}.
\end{array}
\label{e-SCS}
\end{equation}
The algorithm consists of three steps. The first step is projection
onto a subspace, which involves solving a linear system with coefficient
matrix $I+Q$; this is discussed 
in more detail in Section \ref{s-solving-linsys}.
The second step is projection onto a cone, a standard operation discussed
in detail in \cite[\S6.3]{PB:14}. 

The last step is computationally trivial and has a simple interpretation: As
the algorithm runs, the vectors $u^k$ and $\tilde u^k$ converge to each other,
so $u^{k+1}-\tilde u^{k+1}$ can be viewed as the error at iteration $k+1$.  The
last step shows that $v^{k+1}$ is exactly the running sum of the errors.
Roughly speaking, this running sum of errors is used to drive the error to
zero, exactly as in integral control \cite{FPE:94}.

We can also interpret the second and third steps as a combined Moreau
decomposition of the point $\tilde u^{k+1} - v^k$ into its projection onto
$\mathcal C$ (which gives $u^{k+1}$) and its projection onto $-\mathcal C^*$
(which gives $v^{k+1}$).

The algorithm is homogeneous: If we scale the initial points by some factor
$\gamma > 0$, then all subsequent iterates are also scaled by $\gamma$ and the
overall algorithm will give the same primal-dual solution or certificates for
\eqref{e-coneprob}, since the system being solved is also homogeneous.

A straightforward application of ADMM directly to the primal or dual problem in
(\ref{e-coneprob}) obtains an algorithm which requires one linear system solve
involving $A^TA$ and one projection onto the cone $\K$, which has the same
per-iteration cost as (\ref{e-SCS}); see, \eg, \cite{WGY:10} for details.

\subsection{Variations}

There are many variants on the basic ADMM algorithm (\ref{e-SCS})
described above,
and any of them can be employed with the homogeneous self-dual embedding.  We 
briefly describe two important variations that we use in our reference
implementation.

\paragraph{Over-relaxation.}
In the $u$- and $v$-updates, replace all occurrences of $\tilde u^{k+1}$ with 
\[
\alpha \tilde u^{k+1} +(1-\alpha) u^k,
\]
where $\alpha \in {]0,2[}$ is a relaxation parameter \cite{GT:79, EB:92}.
When $\alpha = 1$, this reduces to the basic algorithm given
above.  When $\alpha > 1$, this is known as \emph{over-relaxation};
when $\alpha < 1$, this is \emph{under-relaxation}. Some numerical experiments
suggest that values of $\alpha$ around $1.5$ can improve 
convergence, in practice \cite{Eck:94b, oper_splt_ctrl}.

\paragraph{Approximate projection.}
Another variation replaces the subspace projection update
with a suitable approximation 
\cite{R:76, GT:79, EB:92}.
We replace $\tilde u^{k+1}$ in the first line of (\ref{e-SCS}) 
with any $\tilde u^{k+1}$ that satisfies
\begin{equation}\label{e-approx-proj}
\| \tilde u^{k+1} - (I+Q)^{-1} (u^k + v^k) \|_2 \leq \zeta^k,
\end{equation}
where $\zeta^k>0$ satisfy $\sum_k \zeta^k < \infty$.
This variation is particularly useful when an iterative method is 
used to compute $\tilde u^{k+1}$.

Note that (\ref{e-approx-proj}) is implied by the (more easily
verified) inequality
\begin{equation}\label{e-approx-proj-resid}
\| (Q+I) \tilde u^{k+1} - (u^k + v^k) \|_2 \leq \zeta^k.
\end{equation}
This follows from the fact that $\|(I+Q)^{-1} \|_2 \leq 1$,
which holds since $Q$ is skew-symmetric.
The left-hand side of (\ref{e-approx-proj-resid}) is the norm
of the residual in the equations that define $\tilde u^{k+1}$ 
in the basic algorithm.

\subsection{Convergence}
\label{s-convergence}
\paragraph{Algorithm convergence.}
We show that the algorithm converges, in the sense that it 
eventually produces a point for which the optimality conditions almost 
hold. 
For the basic algorithm (\ref{e-SCS}), and the variant 
with over-relaxation and approximate projection,
for all iterations $k > 0$ we have
\begin{equation}
\label{e-iter-guarantees}
u^k \in \mathcal{C}, \quad v^k \in \mathcal{C}^*, \quad (u^k)^Tv^k = 0.
\end{equation}
These follow from the last two steps of (\ref{e-SCS}),
and hold for any values of $v^{k-1}$ and $\tilde u^k$.
Since $u^{k+1}$ is a projection onto $\mathcal C$, $u^k \in \mathcal C$ 
follows immediately.
The condition $v^k \in \mathcal C^*$ holds since the last
step can be rewritten as $v^{k+1} = \Pi_\mathcal{C^*}(v^k - \tilde u^{k+1})$,
as observed above.
The last condition, $(u^k)^Tv^k =0$, holds by our observation that
these two points are the (orthogonal) Moreau decomposition of the 
same point.

In addition to the three conditions in (\ref{e-iter-guarantees}),
only one more condition must hold for $(u^k,v^k)$ to be optimal:
$Qu^k = v^k$.  This equality constraint holds asymptotically, that is,
we have, as $k \to \infty$,
\begin{equation}\label{e-last-opt-cond}
Qu^k - v^k \to 0.
\end{equation}
(We show this from the convergence result for ADMM below.)
Thus, the iterates $(u^k,v^k)$ satisfy 
three of the four optimality conditions (\ref{e-iter-guarantees}) at every step,
and the fourth one (\ref{e-last-opt-cond}) is satisfied in the limit.

To show that the equality constraint holds asymptotically we use general
ADMM convergence theory; see,
\eg, \cite[\S3.4.3]{BP:11}, or 
\cite{EB:92} for the case of approximate projections.
This convergence theory tells us that 
\begin{equation}\label{e-conv}
\tilde u^k \to u^k,\qquad \tilde v^k \to v^k
\end{equation}
as $k \to \infty$, even with over-relaxation and approximate projection.
From the last step in (\ref{e-SCS}) we conclude that $v^{k+1}-v^k \to 0$.
From~(\ref{e-dsum}), (\ref{e-conv}), and 
$v^{k+1}-v^k \to 0$, we obtain $u^{k+1}-u^k \to 0$.

Expanding (\ref{e-approx-proj-resid}), we have
\[
Q \tilde u^{k+1} + \tilde u^{k+1} - u^k - v^k \to 0,
\]
and using (\ref{e-conv}) we get
\[
Q u^{k+1} + u^{k+1} - u^k - v^k \to 0.
\]
From $u^{k+1}-u^k \to 0$ and $v^{k+1}-v^k \to 0$ we conclude
\[
Qu^k - v^k \to 0,
\]
which is what we wanted to show.

\paragraph{Eliminating convergence to zero.}
We can guarantee that the algorithm will not converge to zero if a
nonzero solution exists, by proper selection of the 
initial point $(u^0, v^0)$,
at least in the case of exact projection.

Denote by $(u^\star, v^\star)$ any nonzero solution to (\ref{e-hsdnice}),
which we assume satisfies either $u^\star_\tau > 0$ or 
$v_\kappa^\star > 0$,
\ie, we can use it to derive an optimal point or a certificate 
for (\ref{e-coneprob}).
If we choose initial point
$(u^0, v^0)$ with $u^0_\tau = 1$ and
$v^0_\kappa = 1$, and all other entries zero, then we have
\[
(u^\star, v^\star)^T (u^0, v^0) > 0.
\]
Let $\phi$ denote the mapping that consists of one iteration of algorithm
(\ref{e-SCS}), \ie, $(u^{k+1}, v^{k+1}) = \phi(u^k, v^k)$. We show in the appendix that the mapping $\phi$ is
nonexpansive, \ie, for any $(u,v)$ and $(\hat u, \hat v)$ we have that 
\begin{equation}
\label{e-nonexpansive}
\|\phi(u, v) - \phi(\hat u, \hat v) \|_2 \leq \| (u ,v) - (\hat u, \hat v) \|_2.
\end{equation}
(Nonexpansivity holds for ADMM more generally; see, \eg, \cite{G:83, FG:83, EB:92} for
details.)
Since $(u^\star, v^\star)$ is a solution to
(\ref{e-hsdnice}) it is a fixed point of $\phi$, \ie,
\begin{equation}
\label{e-phi-fixedpt}
\phi( u^\star, v^\star) = (u^\star, v^\star).
\end{equation}
Since the problem is homogeneous, the point $\gamma (u^\star, v^\star)$ is also a
solution for any positive $\gamma$, and is also a fixed point of $\phi$.
Combining this with (\ref{e-nonexpansive}), we have at iteration $k$ 
\begin{equation}
\label{e-disc}
\|(u^k, v^k) - \gamma (u^\star, v^\star)\|_2^2 \leq \|(u^0, v^0) - \gamma (u^\star, v^\star)\|_2^2,
\end{equation}
for any $\gamma> 0$. Expanding (\ref{e-disc}) and setting
\[
\gamma= \|(u^0,v^0)\|_2^2/ (u^\star, v^\star)^T (u^0, v^0),
\]
which is positive by our choice of $(u^0, v^0)$, we obtain
\[
2(u^\star, v^\star)^T (u^k, v^k) \geq (u^\star, v^\star)^T (u^0, v^0) (1 + \|(u^k, v^k)\|_2^2/\|(u^0, v^0)\|_2^2),
\]
which implies that
\[
(u^\star, v^\star)^T (u^k, v^k) \geq (u^\star, v^\star)^T (u^0, v^0) /2,
\]
and applying Cauchy-Schwarz yields
\begin{equation}
\label{e-nonzero}
\|(u^k, v^k)\|_2 \geq (u^\star, v^\star)^T (u^0, v^0) / 2 \|(u^\star, v^\star)\|_2 > 0.
\end{equation}
Thus, for $k=1,2,\ldots$, the iterates are bounded away from zero.

\paragraph{Normalization.}
The vector given by
\[
(\hat u^k, \hat v^k) = (u^k, v^k) / \|(u^k, v^k)\|_2
\]
satisfies the conditions given in (\ref{e-iter-guarantees}) for all iterations,
and by combining (\ref{e-last-opt-cond}) with (\ref{e-nonzero}) 
we have that 
\[
Q \hat u^k - \hat v^k \rightarrow 0,
\]
in the exact projection case at least.  In other words, the unit vector $(\hat
u^k, \hat v^k)$ eventually satisfies the optimality conditions for the
homogeneous self-dual embedding to any desired accuracy.

\subsection{Termination Criteria}
\label{s-termination}
In view of the discussion of the previous section, a stopping criterion of the
form 
\[
\|Qu^k-v^k\|_2 \leq \epsilon
\]
for some tolerance $\epsilon$, or alternatively a normalized criterion
\[
\|Q u^k - v^k\|_2 \leq \epsilon \|(u^k, v^k)\|_2,
\]
will work, \ie, the algorithm eventually stops.
Here, we propose a different scheme that handles the components of $u$ and $v$
corresponding to primal and dual variables separately. This yields stopping criteria
that are consistent with ones traditionally used for cone programming.

We terminate the algorithm when it finds a primal-dual optimal solution or a certificate
of primal or dual infeasibility, up to some tolerances. If $u_\tau^k>0$,
then let
\[
x^k = u_x^k / u_\tau^k,\quad s^k = v_s^k / u_\tau^k,\quad y^k = u_y^k/u_\tau^k
\]
be the candidate solution. This candidate is guaranteed to satisfy the cone constraints
and complementary slackness condition by \eqref{e-iter-guarantees}. It thus suffices to
check that the residuals
\[
p^k = Ax^{k} + s^{k} - b, \quad
d^k = A^Ty^{k} + c, \quad
g^k = c^T x^k + b^T y^k,
\]
are small. Explicitly, we terminate if
\[
\|p^k\|_2 \leq \epsilon_\mathrm{pri} (1+\|b\|_2), \quad
\|d^k\|_2 \leq \epsilon_\mathrm{dual} (1+\|c\|_2), \quad
|g^k| \leq \epsilon_\mathrm{gap} (1+|c^T x| + |b^T y|)
\]
and emit $(x^k, s^k, y^k)$ as (approximately) primal-dual optimal. Here, quantities
$\epsilon_\mathrm{pri},~\epsilon_\mathrm{dual},~\epsilon_\mathrm{gap}$ 
are the primal residual, dual residual, and
duality gap tolerances, respectively.

On the other hand, if the current iterates satisfy
\[
\|A u_x^k + v_s^k\|_2 \leq (-c^T u_x^k / \|c\|_2) \epsilon_\mathrm{unbdd}, 
\]
then $u_x^k / (-c^T u_x^k)$ is an approximate
certificate of unboundedness with tolerance $\epsilon_\mathrm{unbdd}$, or if they satisfy
\[
\|A^T u_y^k\|_2 \leq (-b^T u_y^k / \|b\|_2) \epsilon_\mathrm{infeas},
\]
then $u_y^k / (-b^Tu_y^k)$ is an approximate
certificate of infeasibility with tolerance $\epsilon_\mathrm{infeas}$. 

These stopping criteria are identical to those used by many other cone solvers
and similar to those used by DIMACS \cite{dimacs, M:03} and the SeDuMi 
solver \cite{sedumi}.

\section{Efficient Subspace Projection}
\label{s-fast-impl}

In this section we discuss how to efficiently compute the projection
onto the subspace $\mathcal Q$, exactly and also approximately (for the
approximate variation).

\subsection{Solving the Linear System}
\label{s-solving-linsys}

The first step is to solve the linear system $(I+Q) \tilde u^k = w$ for some
$w$:
\begin{equation}\label{e-lin-sys}
\begin{bmatrix}
I & A^T & c \\
-A & I & b \\
-c^T & -b^T & 1
\end{bmatrix}
\begin{bmatrix} \tilde u_x\\\tilde u_y\\ \tilde u_\tau \end{bmatrix} = 
\begin{bmatrix} w_x\\w_y\\w_\tau \end{bmatrix}.
\end{equation}
To lighten notation, let
\[
M = \begin{bmatrix} I & A^T \\ -A & I \end{bmatrix}, 
    \quad
h = \begin{bmatrix} c \\ b \end{bmatrix},
\]
so
\[
I + Q = \begin{bmatrix} M & h \\ -h^T & 1 \end{bmatrix}.
\]
It follows that
\[
\begin{bmatrix} \tilde u_x \\ \tilde u_y \end{bmatrix} =
    (M + hh^T)^{-1}\left(\begin{bmatrix} w_x \\ w_y \end{bmatrix} - w_\tau h \right),
\]
where $M + hh^T$ is the Schur complement of the lower right block $1$ in $I + Q$.
Applying the Sherman-Morrison-Woodbury formula \cite[p. 50]{GvL:96} 
to $(M + hh^T)^{-1}$ yields 
\[
\begin{bmatrix} \tilde u_x \\ \tilde u_y \end{bmatrix} =
\left(M^{-1} - \frac{M^{-1} h h^T M^{-1}}{
\left(1+ h^T M^{-1} h \right)}
\right)
\left( \begin{bmatrix} w_x\\w_y \end{bmatrix}  - w_\tau h\right)
\]
and
\[
\tilde u_\tau = w_\tau + c^T \tilde u_x + b^T \tilde u_y.
\]
Thus, in the first iteration, we compute and cache $M^{-1}h$.
To solve \eqref{e-lin-sys} in subsequent iterations, it is only necessary to compute
$M^{-1} (w_x, w_y)$,
which will require the bulk of the computational effort, and then to perform
some simple vector operations using cached quantities. 

There are two main ways to solve linear equations of the form
\begin{equation}
\label{e-linsys}
\begin{bmatrix}
I & -A^T \\
-A & -I  \\
\end{bmatrix}
\begin{bmatrix} z_x \\ -z_y \end{bmatrix}
=
\begin{bmatrix} w_x \\ w_y \end{bmatrix},
\end{equation}
the system that needs to be solved once per iteration.
The first method, a \emph{direct method} that exactly solves the system, is to
solve \eqref{e-linsys} by computing a sparse permuted $LDL^T$ factorization
\cite{davis_book} of the matrix in (\ref{e-linsys})
before the first iteration, then to use this cached
factorization to solve the system in subsequent steps.  This technique, called
factorization caching, is very effective in the common case when the
factorization cost is substantially higher than the subsequent solve cost, so
all iterations after the first one can be carried out quickly.  Because the matrix
is quasi-definite, the factorization is guaranteed to exist for any
symmetric permutation \cite{Van:95}.

The second method, an \emph{indirect method} that we use to approximately
solve the system, involves first rewriting \eqref{e-linsys} as
\[
z_x = (I + A^TA)^{-1}(w_x - A^Tw_y), \quad
z_y = w_y + A z_x,
\]
by elimination. This system is then solved with the conjugate gradient
method (CG) \cite{NW:06,GvL:96,Saad:03}.  Each iteration of conjugate gradient
requires multiplying once by $A$ and once by $A^T$, each of which can be
parallelized.  If $A$ is very sparse, then these multiplications can be
performed especially quickly; when $A$ is dense, it may be better to first form
$G = I + A^T A$ in the setup phase. We warm-start CG by initializing
each subsequent call with the solution obtained by the previous call.
We terminate the CG iterations when the residual satisfies 
(\ref{e-approx-proj-resid}) for some appropriate sequence $\zeta^k$.

\subsection{Repeated Solves}

If the cone problem must be solved more than once, then computation from the
first solve can be re-used in subsequent solves by warm-starting: we set the
initial point to $u^0 = (x^\star, y^\star, 1)$, $v^0 = (0, s^\star, 0)$, where
$x^\star, s^\star, y^\star$ are the optimal primal-dual variables from the
previous solve.  If the data matrix $A$ does not change and a direct method is
being used, then the sparse permuted $LDL^T$ factorization can also be re-used
across solves for additional savings.  This arises in many practical
situations, such as in control, statistics, and sequential convex programming.

\section{Scaling Problem Data}
\label{s-scaling}

Though the algorithm in~\eqref{e-SCS} has no explicit parameters, the relative
scaling of the problem data can greatly affect the convergence.
This suggests a pre-processing step where we scale the data to (hopefully)
improve the convergence.

In particular, consider scaling vectors $b$ and $c$ by positive scalars
$\sigma$ and $\rho$, respectively, and scaling the primal and dual equality
constraints by diagonal positive definite matrices $D$ and $E$,
respectively. This yields the following scaled primal-dual problem pair:
\[
\BA{lcl}
\BA{l}
\BA{ll}
\mbox{\minimize} & \rho (Ec)^T\hat x
\EA \\
\BA{ll}
\mbox{\subjectto} & DA E \hat x + \hat s = \sigma D b \\
& (\hat x, \hat s) \in \reals^n \times \K,
\EA
\EA
& \quad &
\BA{l}
\BA{ll}
\mbox{\maximize} & -\sigma (D b)^T \hat y
\EA \\
\BA{ll}
\mbox{\subjectto} & -EA^T D \hat y +\hat r  = \rho E c \\
& (\hat r, \hat y) \in \{ 0\}^n \times\K^*,
\EA
\EA
\EA
\]
with variables $\hat x$, $\hat y$, $\hat r$, and $\hat s$. After solving this
new cone program with problem data
$\hat A = DAE$, $\hat b = \sigma Db$, and $\hat c = \rho Ec$, the solution
to the original problem~\eqref{e-coneprob} can be recovered
from the scaled solution via
\[
x^\star = E \hat x^\star/\sigma, \quad
s^\star = D^{-1} \hat s^\star / \sigma, \quad
y^\star = D \hat y^\star/\rho.
\]
Transformation by the matrix $D$ must preserve membership of the cone $\K$, to
ensure that if $s \in \K$, then $D^{-1}s \in \K$ (the same is not required of 
$E$). If $\K = \K_1 \times \cdots \times K_q$, where $K_i \in \reals^{m_i}$,
then we could use, for example,
\[
D = \diag(\pi_1 I_{m_1},\ldots, \pi_q I_{m_q}),
\]
where each $\pi_i > 0$.

We have observed that in practice, data which has been \emph{equilibrated},
\ie, scaled to have better conditioning, admits better
convergence \cite{Bauer:63, Bauer:69, Sluis:69, Ruiz:01}. We have found that if
the columns of $A$ and $b$ all have Euclidean norm close to one
and the rows of $A$ and $c$ have similar
norms, then the algorithm~\eqref{e-SCS} typically performs well. The scaling
parameters $E$, $D$, $\sigma$, and $\rho$ can be chosen to (approximately)
achieve this \cite{Osborne:60,Ruiz:01,PC:11}, though the question of whether
there is an optimal scaling remains open.  There has recently been much
work devoted to the question of choosing an optimal, or at least good
diagonal scaling; see \cite{gb:14a, gb:15a}.

\paragraph{Scaled termination criteria.}
When the algorithm is applied to the scaled problem, it is still desirable
to terminate the procedure when the residuals for the \emph{original}
problem satisfy the stopping criteria defined in \S\ref{s-termination}.

The original residuals can be expressed in terms of the scaled
data as 
\[
p^k = (1/\sigma) D^{-1} (\hat A \hat x^{k} + \hat s^{k} - \hat b), \quad
d^k = (1/\rho) E^{-1} (\hat A^T \hat y^{k} + \hat c), \quad
g^k = (1/\rho \sigma)(\hat c^T\hat x^k + \hat b^T \hat y^k),
\]
and the convergence checks can be applied as before. The stopping criteria for
unboundedness and infeasibility then become
\[
\|D^{-1} (\hat A \hat u_x^k + \hat v_s^k )\|_2 \leq (-\hat c^T \hat u_x^k / \|E^{-1} \hat c\|_2) \epsilon_\mathrm{unbdd}, \quad
\|E^{-1} (\hat A^T \hat u_y^k )\|_2 \leq (-\hat b^T \hat u_y^k / \|D^{-1} \hat b\|_2) \epsilon_\mathrm{infeas}.
\]
\section{Numerical Experiments}
\label{s-num-ex}

In this section we present numerical results for SCS, our implementation of the
algorithm described above.  We show results on four application problems, in each
case instances that are small, medium, and large. To demonstrate scaling to
extremely large problems, we also report results on randomly generated problems
with known optimal value.

We compare the results to SDPT3 \cite{SDPT3} and Sedumi \cite{sedumi},
state-of-the-art interior-point solvers.  We use this comparison for several
purposes.  First, the solution computed by these solvers is high accuracy, so
we can use it to assess the quality of the solution found by SCS.  Second, we
can compare the computing times.  Run-time comparison is not completely fair,
since an interior-point method reliably computes a high accuracy solution,
whereas SCS is meant only to compute a solution of modest accuracy and may take
longer than an interior-point method if high accuracy is required. Third,
Sedumi targets the same homogeneous self-dual embedding (\ref{e-hsd}) as SCS,
so we can compare a first-order and a second-order method on the same
embedding.

\subsection{SCS}

Our implementation, which we call SCS for `Splitting Conic Solver', is written
in C and can solve cone programs involving any combination of non-negative,
second-order, semidefinite, exponential, and power cones (and dual exponential
and power cones) \cite{scs}.  It has multi-threaded and single-threaded
versions, and computes the (approximate) projections onto the subspace using
either a direct method or an iterative method. 
SCS is available online at \begin{center}
\url{https://github.com/cvxgrp/scs} \end{center} along with the code to run the
numerical examples. SCS can be used in other C, C++, Python, Matlab, R, Julia,
Java, and Scala programs and is a supported solver in parser-solvers CVX
\cite{cvx}, CVXPY \cite{cvxpy}, Convex.jl \cite{convexjl}, and YALMIP
\cite{yalmip}.  It is now the default solver for CVXPY and Convex.jl for
problems that cannot be expressed using the standard symmetric cones.

The direct implementation uses a single-threaded sparse permuted $LDL^T$
decomposition from the SuiteSparse package \cite{davis_book,ldl,amd}. The
sparse indirect implementation, which uses conjugate gradient, can perform the
matrix multiplications on the CPU or on the GPU. The CPU version uses a basic sparse
multiplication routine parallelized using OpenMP \cite{openmp08}. The
GPU version uses the sparse CUDA BLAS library \cite{cuda}. The indirect solver
uses $\zeta^k = (1/k)^{1.5}$ as the termination tolerance at iteration $k$,
where the tolerance is defined in (\ref{e-approx-proj-resid}).

SCS handles the usual non-negative, second-order, and semidefinite cones,
as well as the exponential cone and its dual
\cite[\S 6.3.4]{PB:14},
\[
K_\mathrm{exp} = \{ (x, y, z) \setsep y > 0,~ ye^{x/y} \leq z \}
\cup \{ (x, y, z) \setsep x \leq 0,~ y = 0,~ z \geq 0 \},
\]
\[
K_\mathrm{exp}^* = \{ (u, v, w) \setsep u < 0,~ -ue^{v/u} \leq ew \}
\cup \{ (0, v, w) \setsep v \geq 0,~ w \geq 0 \},
\]
and the power cone and its dual \cite{nes:06, sy:14, KH:14}, defined as 
\[
K^a_\mathrm{pwr} = \{(x, y, z) \setsep x^a y^{(1-a)} \geq |z|,~ x\geq0,~ y\geq0 \},
\]
\[
(K^a_\mathrm{pwr})^* = \{(u,v,w) \setsep(u/a)^a (v/(1-a))^{(1-a)} \geq |w|,~ u \geq 0,~ v \geq 0\},
\]
for any $a \in [0,1]$.
Projections onto the semidefinite cone are performed using the LAPACK
\verb+dsyevr+ method for computing the eigendecomposition;
projections onto the other cones are implemented in C.
The multi-threaded version computes the projections onto the 
cones in parallel.

In the experiments reported below, 
we use the termination criteria described in \S\ref{s-termination} and
\S\ref{s-scaling}, with the default values
\[
\epsilon_\mathrm{pri} = \epsilon_\mathrm{dual} = \epsilon_\mathrm{gap} =
\epsilon_\mathrm{unbdd} = \epsilon_\mathrm{infeas} = 10^{-3}.
\]
The objective value reported for SCS in the experiments below is the
average of the primal and dual objectives at termination. The time required to
do any preprocessing (such as the matrix factorization) and to carry out
and undo the scaling are included in the total solve times.  

All the experiments were carried out on a system with 32 2.2GHz cores and 512Gb
of RAM, running Linux.  (The single-threaded versions, of course, do not make use
of the multiple cores.) The GPU used was a Geforce GTX Titan X with 12Gb of memory.

\subsection{Lasso}
\label{s-spreg}

Consider the following optimization problem:
\begin{equation}
\label{e-lasso}
\begin{array}{ll}
\mbox{\minimize} & (1/2)\|Fz - g\|_2^2 + \mu \|z\|_1,
\end{array}
\end{equation}
over $z\in \reals^p$, where $F \in \reals^{q\times p}$, $g \in \reals^q$ and
$\mu \in \reals_+$ are data.  This problem, known as
the \emph{lasso} \cite{tibshirani:96}, is widely studied in high-dimensional
statistics, machine learning, and compressed sensing. Roughly speaking,
\eqref{e-lasso} seeks a sparse vector $z$ such that $Fz \approx g$, and the
parameter $\mu$ trades off between quality of fit and sparsity.
It has been observed that first-order methods can perform very well on lasso-type
problems when the solution is sparse \cite{DDD:04, DZ:13}.

The lasso problem can be formulated as the SOCP \cite{lobo:98}
\[
\BA{l}
\begin{array}{ll}
\mbox{\minimize} & (1/2) w + \mu \ones^T t
\EA \\
\BA{ll}
\mbox{\subjectto} & -t \leq z \leq t, \quad 
\left\|\begin{array}{c} 1-w\\
2(Fz - g)
\end{array} \right\|_2 \leq 1 + w
\end{array}
\EA
\]
with variables $z \in \reals^p$, $t \in \reals^p$ and $w \in \reals$.  This
formulation is easily transformed in turn into the standard form
(\ref{e-coneprob}). 

\paragraph{Problem instances.}
We generated data for the numerical instances as follows. First, the entries of
$F$ were sampled independently from a standard normal distribution.  We randomly
generated a sparse vector $\hat z$ with $p$ entries, only $p/10$ of which were
nonzero.  We then set $g = F\hat z + w$, where the entries in $w$ were sampled
independently and identically from $\mathcal{N}(0,0.1)$. We chose $\mu = 0.1
\mu^\mathrm{max}$ for all instances, where $\mu^\mathrm{max} = \|F^Tg\|_\infty$
is the smallest value of $\mu$ for which the solution to (\ref{e-lasso}) is
zero.

\paragraph{Results.}
The results are summarized in Table \ref{t-lasso}.  For the small, medium, and
large instances, the fastest implementation of SCS, indirect on the GPU,
provides a speedup of roughly $30 \times$, $190 \times$, and $1000 \times$,
respectively over SDPT3 and Sedumi.  In the largest case, SCS takes less than
$4$ minutes compared to nearly 3 days for SDPT3 and Sedumi. In other words, not
only is the degree of speedup dramatic in each case, but it also continues to
increase as the problem size gets larger; this is consistent with our goal of
solving problems outside the ability of traditional interior-point methods.

SCS is meant to provide solutions of modest, not high, accuracy.  However, we
see that the solutions returned attain an objective value within 0.01\% of the
optimal value attained by SDPT3 and Sedumi, a negligible difference in
applications. 

If we compare the direct and indirect CPU implementations of SCS, we see that
for small problems the direct version of SCS is faster, but for larger problems
the multi-threaded indirect method dominates. The sparsity pattern in this
problem lends itself to an efficient multi-threaded matrix multiply since the
columns in the data matrix $A$ have a similar number of nonzeros.  This speed-up
is even more pronounced when the matrix multiplications are performed on the GPU.

\begin{table}
\caption{Results for the lasso example.}
\normalsize
\begin{center}
\begin{tabular}{llll}
\hline
& \textbf{small} & \textbf{medium}
& \textbf{large}\\
variables $p$ & 10000 & 30000 & 100000\\
measurements $q$ & 2000 & 6000 & 20000 \\
std.\ form variables $n$ & 2001 & 6001 & 20001\\
std.\ form constraints $m$ & 22002 & 66002 & 220002\\
nonzeros in $A$ & \e{3.8}{6} & \e{3.4}{7} & \e{3.9}{8} \\[1ex]

\textbf{SDPT3}:\\
total solve time& \textbf{196.5 sec} & \textbf{\eb{4.2}{3} sec} & \textbf{\eb{2.3}{5} sec} \\
objective & 682.2 & 2088.0 & 6802.6\\[1ex]

\textbf{Sedumi}:\\
total solve time& \textbf{138.0 sec} & \textbf{\eb{5.6}{3} sec} & \textbf{\eb{2.5}{5} sec} \\
objective & 682.2 & 2088.0 & 6802.6\\[1ex]

\textbf{SCS direct}:\\
total solve time& \textbf{21.9 sec} & \textbf{\eb{3.6}{2} sec} & \textbf{\eb{6.6}{3} sec} \\
factorization time & 5.5 sec & \e{1.1}{2} sec & \e{4.2}{3} sec\\
iterations & 400 & 540 & 500\\
objective & 682.2 & 2088.1 & 6803.5\\[1ex]

\textbf{SCS indirect}:\\
total solve time & \textbf{31.6 sec} & \textbf{\eb{1.2}{2} sec} & \textbf{\eb{7.5}{2} sec} \\
average CG iterations & 5.9 & 5.9 & 5.9\\
iterations & 400 & 540 & 500 \\
objective & 682.2 & 2088.1 & 6803.6 \\[1ex]

\textbf{SCS indirect GPU}:\\
total solve time & \textbf{4.6 sec} & \textbf{22.0 sec} & \textbf{\eb{2.1}{2} sec} \\
\hline
\end{tabular}
\end{center}
\label{t-lasso}
\end{table}

\subsection{Portfolio Optimization}

Consider a simple long-only portfolio optimization problem
\cite{Markowitz1952,port_opt_bound}, \cite[\S 4.4.1]{BV:04}, in which we choose
the relative weights of assets to maximize the expected risk-adjusted return of a
portfolio:
\[
\BA{ccccc}
\mbox{\maximize} & \left[ \mu^Tz - \gamma(z^T \Sigma z) \right] & & \mbox{\subjectto} & \ones^Tz = 1 , \quad z \geq 0,
\EA
\]
where the variable $z\in \reals^p$ represents the portfolio of $p$
assets,
$\mu \in \reals^p$ is the vector of expected returns, $\gamma>0$ is
the \emph{risk aversion parameter}, and
$\Sigma \in \reals^{p \times p}$ is the asset return covariance matrix
(also known as the \emph{risk model}).
The risk model is expressed in \emph{factor model form}
\[
\Sigma = FF^T + D,
\]
where $F \in \reals^{p\times q}$ is the \emph{factor loading matrix}
and $D \in \reals^{p\times p}$ is a diagonal
matrix representing `idiosyncratic' or asset-specific risk.
The number of risk factors $q$ is typically much less
than the number of assets $p$. 
(The factor model form is widely used in practice.)

This problem can be converted in the standard way into an SOCP:
\begin{equation}
\label{e-port}
\BA{l}
\begin{array}{ll}
\mbox{\maximize} & \mu^Tz - \gamma(t + s) 
\EA \\
\BA{ll}
\mbox{\subjectto} & \ones^Tz =1, \quad z \geq 0, \quad \| D^{1/2} z \|_2 \leq u, \quad \| F^T z \|_2 \leq v \\
& \left\|\left(1 - t, 2u\right)\right\|_2 \leq 1+t, \quad
\left\|\left(1 - s, 2v\right)\right\|_2 \leq 1+s,
\end{array}
\EA
\end{equation}
with variables $z \in \reals^p$, $t \in \reals$, $s \in \reals$, $u \in
\reals$, and $v \in \reals$. This can be transformed into standard form
\eqref{e-coneprob} in turn.

\paragraph{Problem instances.}
The vector of log-returns, $\log(\mu)$, was sampled from a standard normal
distribution, yielding log-normally distributed returns.  The entries in $F$
were sampled independently from $\mathcal{N}(0,0.1)$, and the diagonal entries of $D$ were sampled
independently from a uniform distribution on $[0,0.1]$.  For all problems, we
chose $\gamma = 1$. 

\paragraph{Results.}
The results are summarized in Table \ref{t-port}.
In all cases the objective value
attained by SCS was within $0.5\%$ of the optimal value. The worst budget
constraint violation of the solution returned by SCS in any instance was only
$0.002$ and the worst non-negativity constraint violation was only
\e{5}{-7}. SCS direct is more than $7$ times faster than SDPT3 on the
largest instance, and much faster than Sedumi, which didn't manage to solve
the largest instance after a week of computation.

Unlike the previous example, the direct solver is faster than the indirect
solver on the CPU for all instances. This is due to imbalance in the number of
nonzeros per column which, for the simple multi-threaded matrix multiply we're
using, leads to some threads handling much more data than others, and so the
speedup provided by parallelization is modest. The indirect method on the GPU
is fastest for the medium sized example. For the small example the cost of
transferring the data to the GPU outweighs the benefits of performing the
computation on the GPU, and the large example could not fit into the GPU memory.

\begin{table}
\caption{Results for the portfolio optimization example.}
\normalsize
\begin{center}
\begin{tabular}{llll}
\hline
 & \textbf{small} & \textbf{medium} 
& \textbf{large}\\
assets $p$ & 100000 & 500000 & 2500000 \\
factors $q$ & 100 & 500 & 2500 \\
std.\ form variables $n$ & 100103 & 500503 & 2502503 \\
std.\ form constraints $m$ & 200104 & 1000504 & 5002504 \\
nonzeros in $A$ & \e{1.3}{6} & \e{2.5}{7} & \e{5.1}{8} \\[1ex]

\textbf{SDPT3}:\\
total solve time & \textbf{70.7 sec} & \textbf{\eb{1.6}{3} sec} & \textbf{\eb{6.3}{4} sec} \\
objective & 0.0388 & 0.0364 & 0.0369 \\[1ex]

\textbf{Sedumi}:\\
total solve time & \textbf{100.6 sec} & \textbf{\eb{7.9}{3} sec} & \textbf{$>$ \eb{6.1}{5} sec} \\
objective & 0.0388 & 0.0364 & ? \\[1ex]

\textbf{SCS direct}:\\
total solve time& \textbf{13.0 sec} & \textbf{190 sec} & \textbf{\eb{9.6}{3} sec} \\
factorization time & 0.6 sec & 19.2 sec & 913 sec\\
iterations & 500 & 440 & 980 \\
objective & 0.0388 & 0.0365 & 0.0367 \\[1ex]

\textbf{SCS indirect}:\\
total solve time & \textbf{27.6 sec} & \textbf{313 sec} & \textbf{\eb{2.5}{4} sec} \\
average CG iterations & 3.0 & 3.0 & 3.0 \\
iterations & 500 & 440 & 980 \\
objective & 0.0388 & 0.0365 & 0.0367 \\
\textbf{SCS indirect GPU}:\\
total solve time & \textbf{27.8 sec} & \textbf{184 sec} & \textbf{OOM} \\
\hline
\end{tabular}
\end{center}
\label{t-port}
\end{table}

\subsection{Robust Principal Components Analysis}

This example considers the problem of recovering a low rank matrix
from measurements that have been corrupted by sparse noise \cite{rpca, AGM:14}.
In \cite{rpca}, the authors formulated this problem as follows:
\begin{equation}
\label{e-rpca}
\begin{array}{ll}
\mbox{\minimize} & \|L\|_* \\
\mbox{\subjectto} & \|S\|_1 \leq \mu \\
& L + S = M,
\end{array}
\end{equation}
with variables $L \in \reals^{p\times q}$ and $S \in \reals^{p \times q}$, and with data
$M \in \reals^{p \times q}$ the matrix of measurements and 
$\mu \in \reals_+$
a parameter that constrains the estimate of the corrupting noise term
to be below a certain value.  Here, $\|\cdot\|_*$ is the nuclear norm (dual of spectral norm) 
and $\| \cdot \|_1$ is the elementwise $\ell_1$ norm (\ie, sum of the 
absolute values of the entries).  Roughly
speaking, the problem is a convex surrogate for decomposing the 
given matrix $M$ into the sum of a
sparse matrix $S$ and a low-rank matrix $L$.

The problem can be converted into an SDP as follows \cite{FHB:01, VB:96}:
\[
\begin{array}{ll}
\mbox{\maximize} & (1/2)(\Tr(W_1) + \Tr(W_2))\\
\mbox{\subjectto} & -t \leq {\bf vec}(S) \leq t \\
& L + S = M \\
& \ones^T t \leq \mu \\
& \left[
\begin{array}{cc}
W_1 & L \\ L^T & W_2
\end{array}
\right] \succeq 0
\end{array}
\]
with variables $L \in \reals^{p \times q}$, $S \in \reals^{p \times q}$, 
$W_1 \in \reals^{p \times p}$, $W_2 \in \reals^{q \times q}$, and $t \in \reals^{pq}$,
where ${\bf vec}(S)$ returns the columns of $S$ stacked as a single vector.
The transformation of this problem into standard form (\ref{e-coneprob})
is straightforward.

\paragraph{Problem instances.}
We set $M = \hat L + \hat S$ where $\hat L$ was a randomly generated rank-$r$
matrix and $\hat S$ was a sparse matrix with approximately $10\%$ nonzero entries.
For all instances, we set $\mu$ to be equal to the sum of absolute values of
the entries of $\hat S$ and generated the data with $r = 10$. For simplicity, we
chose the matrices to be square, \ie, $p=q$, for all instances.

\paragraph{Results.}
The results are summarized in Table \ref{t-rpca}.
On the two larger examples, SDPT3 and Sedumi both ran out of memory, so we
cannot directly measure the suboptimality of the SCS solution. However, the
reconstruction error
\[
\|L - \hat L\|_* / \|\hat L \|_*,
\]
where $\hat L$ is the true low-rank matrix used to generate the data and $L$ is
the estimate returned by our algorithm, was less than $3 \times 10^{-4}$ across all
instances. Since this is the actual metric of interest in applications, this
implies that the solutions returned were more than adequate.

In this example the direct, indirect, and indirect GPU implementations of SCS
take roughly the same amount of time. This is because the time required to
project onto the semidefinite cone is the dominant cost per iteration 
(for the medium and large problems), rather than the linear system solve.  


\begin{table}
\caption{Results for the robust PCA example.}
\normalsize
\begin{center}
\begin{tabular}{llll}
\hline
 & \textbf{small} & \textbf{medium} 
& \textbf{large}\\
matrix dimension $p$ & 100 & 500 & 1000\\
std.\ form variables $n$ & 10001 & 250001  & 1000001 \\
std.\ form constraints $m$ & 40101 & 1000501 & 4001001 \\
nonzeros in $A$ & $5.0 \times 10^4$ & $1.3 \times 10^6$ & $5.0 \times 10^6$ \\[1ex]

\textbf{SDPT3}:\\
total solve time & \textbf{429.1 sec} & \textbf{OOM} & \textbf{OOM} \\
objective & 959.3 & OOM & OOM \\[1ex]

\textbf{Sedumi}:\\
total solve time & \textbf{\eb{9.0}{3} sec} & \textbf{OOM} & \textbf{OOM} \\
objective & 959.3 & OOM & OOM \\[1ex]

\textbf{SCS direct}:\\
total solve time& \textbf{2.0 sec} & \textbf{18.5 sec} & \textbf{94.8 sec} \\
factorization time & \e{8.1}{-2} sec & \e{9.9}{-1} sec & 3.7 sec \\
iterations & 120 & 60 & 80 \\
objective & 959.3 & $4.9 \times 10^3$ & $1.0 \times 10^4$ \\[1ex]

\textbf{SCS indirect}:\\
total solve time & \textbf{1.8 sec} & \textbf{16.4 sec} & \textbf{92.3 sec} \\
average CG iterations & 1.0 & 1.0 & 1.0\\
iterations& 120 & 60 & 80\\
objective & 959.3 & $4.9 \times 10^3$ & $1.0 \times 10^4$ \\
\textbf{SCS indirect GPU}:\\
total solve time & \textbf{6.1 sec} & \textbf{18.2 sec} & \textbf{125.0 sec} \\
\hline
\end{tabular}
\end{center}
\label{t-rpca}
\end{table}

\subsection{Logistic Regression with $\ell_1$-Regularization}

In logistic regression the goal is to find the maximum likelihood fit of a
logistic model to (binary) labeled data 
\cite[\S7.1.1]{BV:04}. In this problem we add an additional
regularization term, which increases
the sparsity of the solution. A fixed parameter $\mu \geq 0$ trades off the
likelihood of the model and the model sparsity.

Given data points $z_1, \ldots, z_q \in \reals^p$, with labels $y_1, \ldots,
y_q \in \{-1,1\}$, the $\ell_1$-regularized logistic regression problem is
given by \cite{S:10,FHT:10}
\begin{equation}
\label{e-logregorig}
\begin{array}{ll}
\mbox{\minimize} & \sum_{i=1}^q \log(1 + \exp(y_i w^T z_i)) + \mu \|w \|_1 \\
\end{array}
\end{equation}
with variable $w \in \reals^p$. 

This problem can be converted into a convex cone problem over 
a product of exponential cones
$K_\mathrm{exp} \subset \reals^3$ as follows 
\begin{equation}
\label{e-logreg}
\begin{array}{ll}
\mbox{\minimize} & \ones^T t + \mu \ones^T s \\
\mbox{\subjectto} & -s \leq w \leq s \\ 
& u + v \leq 1 \\
&\begin{bmatrix} y_iw^Tz_i - t_i \\ 1 \\ u_i
\end{bmatrix} \in K_\mathrm{exp},\quad i =1, \ldots, q \\[15pt]
&\begin{bmatrix} -t_i \\ 1 \\ v_i
\end{bmatrix} \in K_\mathrm{exp},\quad i =1, \ldots, q,
\end{array}
\end{equation}
with variables $w \in \reals^p$, $s \in \reals^p$, $t \in \reals^q$, $u \in \reals^q$, and $v \in \reals^q$,
which is readily transformed into standard form (\ref{e-coneprob}).

\paragraph{Problem instances.}
The data were generated as follows. First, we randomly selected a weights
vector $w_\mathrm{true} \in \reals^p$ with at most $p/5$ of the entries
nonzero. Then, each data point $z_i$ was sampled from a standard normal
distribution and assigned a positive label with probability
equal to the value of the logistic function applied to $w_\mathrm{true}^Tz_i$.
For each instance we set $\mu = 0.1 \mu^\mathrm{max}$, where $\mu^\mathrm{max}
= (1/2) \|\sum_{i=1}^q y_i z_i \|_\infty$ is the smallest value of $\mu$ for
which the solution to (\ref{e-logregorig}) is zero.  

\paragraph{Results.}
The results are summarized in table \ref{t-logreg}.  Neither SDPT3 nor Sedumi
can solve exponential cone programs, so we cannot make a direct comparison
between SCS and the interior point solvers in this case. However, using CVX we
can approximate an exponential cone program using a sequence of SDPs. With this
technique SDPT3 is able to solve the smallest instance in a little under two
hours, achieving an objective value of $3876.97$, a difference of less than
$0.001\%$ when compared to SCS on the same problem. Despite this, SCS is able
to solve the largest instance, with almost a billion nonzeros in the data
matrix, in just a few hours using both direct and indirect solvers.

In this example the sparsity pattern of the data matrix does not lend itself to
efficient multi-threaded matrix multiplies when parallelizing over columns.
Because of this the indirect method has little advantage over the direct method.
The indirect method on the GPU is the fastest solver for the small and medium sized
problems, but the GPU did not have enough memory to solve the large instance.

\begin{table}
\caption{Results for the $\ell_1$-regularized logistic regression example.}
\normalsize
\begin{center}
\begin{tabular}{llll}
\hline
 & \textbf{small} & \textbf{medium} 
& \textbf{large}\\
features $p$ & 100 & 1000 & 10000 \\
samples $q$ & 10000 & 100000 & 1000000 \\
std.\ form variables $n$ &30200 & 302000 & 3020000 \\
std.\ form constraints $m$ & 70200 & 702000 & 7020000 \\
nonzeros in $A$ & $1.6 \times 10^5$ & $1.0 \times 10^7$ & $9.2 \times 10^8$\\[1ex]

\textbf{SCS direct}:\\
total solve time& \textbf{5.4 sec} & \textbf{824 sec} & \textbf{\eb{2.5}{4} sec} \\
factorization time & 0.7 sec & 683 sec & \e{5.9}{3} sec\\
iterations & 280 & 380 & 860 \\
objective & 3877.0  & \e{1.2}{4} & \e{4.3}{4} \\[1ex]

\textbf{SCS indirect}:\\
total solve time & \textbf{6.2 sec} & \textbf{230 sec} & \textbf{\eb{3.7}{4} sec} \\
average CG iterations & 3.29 & 3.39 & 4.06 \\
iterations& 280 & 380 & 860 \\
objective & 3877.1 & \e{1.2}{4} & \e{4.3}{4} \\

\textbf{SCS indirect GPU}:\\
total solve time & \textbf{4.4 sec} & \textbf{120 sec} & \textbf{OOM} \\
\hline
\end{tabular}
\end{center}
\label{t-logreg}
\end{table}

\subsection{Random Cone Programs}
\label{s-random-cone}
In this subsection we describe how to generate a feasible bounded 
random cone program
with known optimal objective value, given problem dimensions $n$ and $m$ and
cone $\K$, and present SCS performance results on three random SOCPs. The
procedure simultaneously generates the data $(A,b,c)$ and a primal-dual
solution $(x^\star, s^\star, y^\star)$.
The solution need not be unique, so we do not expect to 
recover $(x^\star, s^\star, y^\star)$; we do expect to recover nearly 
primal and dual feasible 
points with nearly the same objective value.

First, we generate a random vector $z \in \reals^m$ and set $s^\star =
\Pi_\K(z)$ and $y^\star = s^\star - z$.  This ensures conic feasibility,
complementary slackness, and a zero duality gap by Moreau.
Next we randomly generate the data matrix $A\in \reals^{m \times n}$, 
with any desired sparsity pattern, and randomly generate the primal
solution $x^\star \in \reals^n$. Finally, we set $b = Ax^\star + s^\star$ and $c
= -A^T y^\star$, which ensures equality constraint feasibility. The solution to
the problem is not necessarily unique, but the optimal value is given by
$c^T x^\star$.

A similar procedure can be used to generate infeasible or unbounded random cone
programs by simultaneously generating the problem data and a certificate of
primal or dual infeasibility. However, it is not as easy to ensure that the data
matrix is sparse in those cases.

We use this method to generate three SOCPs of different sizes.
The data matrix $A$ was
generated by selecting the nonzero entries uniformly at random, and generating
the nonzero values by sampling from a standard normal distribution.  
Even the small instance is large; the large instance involves more 
than 100Gb of data, and is extremely large.
We used the indirect linear system solver to solve these three problem
instances.

\paragraph{Results.}
The results are given in table \ref{t-random-socps}. The
results indicate that even very large problems can be solved to modest accuracy
with just a few thousand applications of the data matrix and its adjoint.

\begin{table}
\caption{Results for randomly generated cone programs.}
\normalsize
\begin{center}
\begin{tabular}{llll}
\hline
& \textbf{small} & \textbf{medium} & \textbf{large}\\
variables $n$ & \e{1}{4} & \e{1}{6} & \e{4.8}{6}\\
constraints $m$ & \e{3}{4} & \e{3}{6} & \e{1.4}{7}\\
nonzeros in $A$ & \e{1}{6} & \e{1}{9} & \e{1.1}{10}\\
size of $A$ & 11.2Mb & 11.2Gb & 156.8Gb\\[1ex]
\textbf{SCS indirect}:\\
total solve time & 0.9 sec & \e{4.3}{3} sec & \e{1.9}{5} sec\\
iterations &40 & 160 & 240\\
average CG iterations &4.7 & 6.1 & 6.1\\
total matrix multiplies & 457 & 2265 & 3408\\
$|c^T x - p^\star|/|p^\star|$ & \e{2.0}{-3} & \e{2.0}{-3} & \e{1.8}{-3} \\
$|b^T y - p^\star|/|p^\star|$ & \e{1.2}{-4}& \e{1.3}{-4} & \e{3.4}{-4} \\
\end{tabular}
\end{center}
\label{t-random-socps}
\end{table}

\section{Conclusions}

We presented an algorithm that can return primal and dual optimal points for
convex cone programs when possible, and certificates of primal or dual
infeasibility otherwise. The technique involves applying an
operator splitting method, the alternating direction method of multipliers, to
the homogeneous self-dual embedding of the original optimization problem.  This
embedding is a feasibility problem that involves finding a point in the
intersection of an affine set and a convex cone, and each iteration of our
method solves a system of linear equations and projects a point onto the cone.
We showed how these individual steps can be implemented efficiently and are often
amenable to parallelization. 
We discuss methods for automatic problem scaling, 
a critical step in making the method robust.

We provide a reference implementation of our algorithm in C, which we
call SCS. We show that this solver can solve large instances of cone problems
to modest accuracy quickly and is particularly well suited to solving large
cone problems outside of the reach of standard interior-point methods.
As far as we know, the problems reported in \S\ref{s-random-cone}
are the largest general purpose cone problems solved to date.

\begin{acknowledgements}
This research was supported by DARPA's XDATA program under grant
FA8750-12-2-0306. N. Parikh was supported by a NSF Graduate Research Fellowship
under grant DGE-0645962. 
The authors thank Wotao Yin for extensive comments and suggestions on an
earlier version of this manuscript, and Lieven Vandenberghe for fruitful
discussions early on. We would also like to thank the anonymous reviewers for
their constructive feedback.
\end{acknowledgements}

\appendix
\setcounter{secnumdepth}{0}
\section{Appendix: Nonexpansivity}
\label{s-app-nonexp}
In this appendix we show that the mapping consisting of one iteration
of the algorithm (\ref{e-SCS}) is nonexpansive, \ie, if we denote
the mapping by $\phi$, then we shall show that
\[
\|\phi(u, v) - \phi(\hat u, \hat v) \|_2 \leq \| (u ,v) - (\hat u, \hat v) \|_2,
\]
for any $(u,v)$ and $(\hat u, \hat v)$. 

From (\ref{e-SCS}) we can write the
mapping as the composition of two operators, $\phi = P \circ L$, where
\[
P(x) = (\Pi_\mathcal{C}(x), -\Pi_{-\mathcal{C}^*}(x)),
\]
and 
\[
L(u,v) = (I+Q)^{-1}(u+v) - v.
\]
To show that $\phi$ is nonexpansive we only need to show that
both $P$ and $L$ are nonexpansive.

To show that $P$ is nonexpansive we proceed as follows
\begin{align*}
\|x - \hat x\|_2^2 &= \|\Pi_\mathcal{C}(x) + \Pi_{-\mathcal{C}^*}(x) - \Pi_\mathcal{C}(\hat x) - \Pi_{-\mathcal{C}^*}(\hat x) \|_2^2 \\ 
&= \|\Pi_\mathcal{C}(x)- \Pi_\mathcal{C}(\hat x)\|_2^2 + \|\Pi_\mathcal{-C^*}(x)- \Pi_\mathcal{-C^*}(\hat x)\|_2^2 \\ 
& \qquad-2 \Pi_\mathcal{C}(\hat x)^T \Pi_\mathcal{-C^*}(x) - 2 \Pi_\mathcal{C}(x)^T \Pi_\mathcal{-C^*}(\hat x) \\
&\geq \|\Pi_\mathcal{C}(x)- \Pi_\mathcal{C}(\hat x)\|_2^2 + \|\Pi_\mathcal{-C^*}(x)- \Pi_\mathcal{-C^*}(\hat x)\|_2^2 \\
&=  \| (\Pi_\mathcal{C}(x)- \Pi_\mathcal{C}(\hat x) ) , - (\Pi_\mathcal{-C^*}(x)- \Pi_\mathcal{-C^*}(\hat x) )\|_2^2 \\
&= \|P(x) - P(\hat x) \|_2^2,
\end{align*}
where the first equality is from the Moreau decompositions of $x$ and $\hat x$
with respect to the cone $\mathcal{C}$, the
second follows by expanding the norm squared and the fact that
$\Pi_\mathcal{C}(x) \perp \Pi_\mathcal{-C^*}(x)$ for any $x$, and the inequality
follows from $\Pi_\mathcal{C}(\hat x)^T \Pi_\mathcal{-C^*}(x) \leq 0$ by the
definition of dual cones. 

Similarly for $L$ we have
\begin{align*}
\|L(u,v) - L(\hat u, \hat v) \|_2  &=  \left\|(I+Q)^{-1}(u - \hat u +v - \hat v) - v + \hat v\right\|_2 \\
&= \left\| \begin{bmatrix}  (I+Q)^{-1} & -(I-(I+Q)^{-1}) \end{bmatrix} (u - \hat u, v - \hat v) \right\|_2 \\
&\leq \|(u - \hat u, v - \hat v) \|_2 = \|(u,v) - (\hat u, \hat v)\|_2,
\end{align*}
where the inequality can be seen from the fact that 
\[
\begin{bmatrix}  (I+Q)^{-1} & -(I-(I+Q)^{-1}) \end{bmatrix} \begin{bmatrix}  (I+Q)^{-1} & -(I-(I+Q)^{-1}) \end{bmatrix}^T = I
\]
by the skew symmetry of $Q$, and so 
$\left\|\begin{bmatrix}  (I+Q)^{-1} & -(I-(I+Q)^{-1}) \end{bmatrix}\right\|_2 = 1$.

\paragraph{Conflict of interest.} 
The authors declare that they have no conflict of interest.

\newpage

\bibliographystyle{spmpsci_unsrt}
\bibliography{scs}

\end{document}